

\newif\ifprima \primatrue \newtoks\piede \piede={\hfill}
\newtoks\altre \altre={\hfill\tenrm\folio\hfill}
\footline={\ifprima\the\piede\global\primafalse\else\the\altre\fi} 

\hsize=125mm\vsize=187mm\hoffset=0mm\voffset=3mm
\overfullrule=0pt\pretolerance=500\tolerance=1000\brokenpenalty=5000
\magnification\magstep 1 

\font\tit=cmbx12
\font\ma=cmcsc10
\font\mi=cmcsc8
\font\otrm=cmr8 \font\seirm=cmr6
\font\otb=cmbx8
\font\otsi=cmsy8 \font\seisi=cmsy6
\font\otmi=cmmi8 \font\seimi=cmmi6
\def\ottopunti{%
\textfont0=\otrm \scriptfont0=\seirm%
\def\rm{\fam0\otrm}%
\textfont1=\otmi \scriptfont1=\seimi%
\textfont2=\otsi \scriptfont2=\seisi%
\textfont\bffam=\otb%
\abovedisplayskip=9pt plus 2pt minus 6pt
\belowdisplayskip=\abovedisplayskip
\abovedisplayshortskip=0pt plus 2pt
\belowdisplayshortskip=5pt plus 2pt minus 3pt
\smallskipamount=2pt plus 1pt minus 1pt
\medskipamount=4pt plus 2pt minus 2pt
\bigskipamount=9pt plus 4pt minus 4pt
\normalbaselineskip=9pt
\normalbaselines 
\rm}


\def\bib{\hangindent=35pt\hangafter=1}


\font\seirm=cmr6
\font\seisi=cmsy6
\font\seimi=cmmi6


\def\fine{\unskip\kern 6pt\penalty 500\raise -1pt\hbox{\vrule\vbox 
to8pt{\hrule width 6pt \vfill\hrule}\vrule}\par}

\def\fr{\rightarrow}
\def\bfr{\longrightarrow}
\def\Hc#1#2#3#4{H^{#1}({\cal{#2}}_{#3}({#4}))}
\def\h#1#2#3#4{h^{#1}({\cal{#2}}_{#3}({#4}))}

\def\H*#1#2#3{H^{#1}_*({\cal{#2}}_{#3})} 
\def\Ho#1#2#3{H^{#1}({\omega}_{#2}({#3}))}

\def \Hom {{\rm Hom\,}}
\def \Ext {{\rm Ext\,}}

\def\mapright#1{\smash{\mathop{\longrightarrow}\limits^{#1}}}
\def\grf#1#2{\overbrace{{\hskip#1cm}}^{#2}}

\def \remk#1 {\medbreak \noindent {\bf Remark #1.}\enspace }
\def \remks#1 {\medbreak \noindent {\bf Remarks #1.}\enspace }
\def \exp#1 {\medbreak \noindent {\bf Example #1.}\enspace }

\def \bP {{\bf P}}
\def \cB {{\cal B}}
\def \cA {{\cal A}}

\def \cO {{\cal O}}

\def \cJ {{\cal J}}
\def \cM {{\cal M}}
\def \cN {{\cal N}}

\def \bZ {{\bf Z}}


\null

\vskip 1in
\centerline{\tit Biliaison classes of curves in ${\bf P}^3$}
\vskip 1cm
\centerline{{\ma Rosario Strano}  \footnote {(*)}{1991 {\it
Mathematics Subject Classification:} primary 14H50, secondary 14H45.}
\footnote{}{Work done with the
financial support of the Ministry of  Scientific Research.}}
\vskip 1cm
\noindent


{\noindent\narrower
{\mi Abstract}.\thinspace {\ottopunti
We characterize the curves in ${\bf P}^3$ which are minimal in their biliaison class. Such curves are exactly the curves which do no admit an elementary descending biliaison. As a consequence we have that every curve in ${\bf P}^3$ can be obtained from a minimal one by means of a finite sequence of ascending elementary biliaisons. 
}
\par}

\bigskip\noindent


\bigskip
\noindent
{\bf 0. Introduction.}  Let $C$ be a curve ( i.e. a locally C.M.
subscheme equidimensional of dimension 1) in
${\bf P}={\bf P}^3_{\bf k}$ over an algebraically closed field $\bf k$ . Let $\cJ_C$  be its ideal sheaf in $\cO_\bP$ and $I_C=\H*0JC$ be the homogeneous 
ideal of $C$  in the polynomial ring $R={\bf 
k}[x_0,\dots,x_3]=\H*0O{\bf P}$. Denote by $$ s(C)={\rm inf }\{n\in \bZ\ |\ \h0JCn \neq 0\};$$
$$ e(C)={\rm sup }\{n\in \bZ\ |\ \h1OCn \neq 0\}.$$
In [R] P. Rao introduced the notion of biliaison (double linkage) class of curves and proved that two curves $C$ and $C'$ are in the same classe if and only if $\H*1JC \cong \H*1J{C'}(h)$ for some $h\in \bZ$. In [LR] R. Lazarsfeld and P. Rao proved that curves with $s(C)\geq e(C)+4$ are minimal in their biliaison class. In [MDP] M. Martin-Deschamps and D. Perrin gave a costruction of the minimal curve in each biliaison class starting from the minimal free graded resolution of the Rao module $M(C)=\H*1JC$.\par
In this paper, by completing the result of Lazarsfeld and Rao, we give a characterization of minimal non ACM curves in the case $s(C)\leq e(C)+3$; these curves are those satisfying the following condition:
(*) for every  $s(C)\leq s \leq e(C)+3$ there is a form $H_s$ of degree $<s$ dividing every 
element in $\Hc0JCs$ and such that $H_s\cdot\Ho0C{3-s}=0$.\par
We see also that these curves are exactly the curves that do not admit any elementary descending biliaison; from this one can deduce the following theorem that improves
a result obtained by M. Martin-Deschamps and Perrin [ MDP] and by E. Ballico, G. Bolondi and J. Migliore [BBM]:\par
{\it Every curve in ${\bf P}^3$ can be obtained from a minimal one by means of a finite number of ascending elementary biliaisons}.\par
The paper is organized as follows: in Sections 1 and 2 we prove the main result (Theorem 1) and precisely in Section 1 we prove that a non ACM curve with $s(C)\leq e(C)+3$ that do not admit any elementary descending biliaison satisfy condition (*), while in Section 2  we prove that such curve is minimal.\par
In Section 3 we deduce the above result regarding the structure of a given biliaison class together with some other interesting consequences.\par
 We refer to [MDP] for all the notations and in particular for the basic results  concerning the biliaison.\par
By abuse of notation we will denote with the same symbol an homogeneous polinomial in $R$ and the surface of ${\bf P}^3$ that it defines.

\par
\bigskip

\noindent
{\bf 1.} In this section first we state the main result of the paper. 
 The notations are as in the introduction. 
In order to include all possible curves in ${\bf P}^3$ we consider a line
as  a minimal curve in the biliaison class of the ACM curves in ${\bf P}^3$.
\proclaim{\bf Theorem 1}. 
 For a curve $C\subset {\bf P}^3$ the following are equivalent:
{\item 1\kern -5pt)} $C$ is minimal in its biliaison class
{\item 2\kern -5pt)} $C$ does not admit any elementary descending biliaison
{\item 3\kern -5pt)}  $C$ is of one (and only one) of the following types:
{\itemitem a\kern -5pt)} a line
{\itemitem b\kern -5pt)} $s(C)\geq e(C)+4$ 
{\itemitem c\kern -5pt)} $s(C)\leq e(C)+3 $  and
{\itemitem \rm{(*)}} for every  $s(C)\leq s \leq e(C)+3$ there is a form $H_s$ of positive degree $<s$ dividing every 
element in $\Hc0JCs$ and such that $H_s\cdot\Ho0C{3-s}=0$.\par\bigskip\noindent
{\it Proof}. In the rest of this Section we prove that 2) implies 3). First we prove some preliminary result; recall that, for every $Q\in\Hc0JCs$ and for every $m\in\bZ$ we have an exact sequence:

$$ 0\fr \Hc0OQ{m}\fr \Hom (\cJ_{C/Q}, \cO_Q(m)) \fr \Ho0C{4-s+m}\fr 0 $$\noindent

in particular for $m=-1$ we have an isomorphism : $\Hom (\cJ_{C/Q}, \cO_Q(-1))\cong\Ho0C{3-s}$. 

\proclaim{\bf Proposition 1}. Let $s(C)\leq s \leq e(C)+3$ ,
$Q\in\Hc0JCs$,  $\xi\in \Ho0C{3-s}$, $\xi\neq 0$ and let $\eta: \cJ_{C/Q}\fr \cO_Q(-1)$ the morphism of 
$\cO_Q$-modules corresponding to $\xi$. Then $\eta$ is not injective if and only if there is a form $H$ of degree $<s$ dividing $Q$ and such that $H\cdot\xi=0$.\par
{\it Proof of Proposition 1}. The proof is similar to that  given in [S] Proposition 2.2.3; see also [MDP] Proposition III.2.6.

 Assume that  $\eta$ is not injective; by [S] Proposition 2.2.3 its kernel and its image are respectively of the form $\cJ_{D/Q}$ and $\cJ_{E/Q}(-1)$, where $D$,$E$ are subschemes of $Q$ containing respectively surfaces $H$ , $ K$ of degree $<s$ such that $H\cdot K=Q$. Since  $\cJ_{E/Q}(-1)\subseteq\cJ_{K/Q}(-1)$ and $H$ kills $\cJ_{K/Q}(-1)$ we have $H\cdot \eta=0$ and hence  $H\cdot\xi =0$. 
Conversely assume that there is a form $H$ of degree $<s$ dividing $Q$ and such that $H\cdot\xi=0$ and let denote by 
$u: I_{C/Q}\fr R_Q(-1)$ the  (degree zero) homomorphism of graded $R$-modules corresponding to $\eta$; 
 from the exact sequence 
$$ 0\fr \Hc0OQ{h-1}\fr \Hom (\cJ_{C/Q}, \cO_Q(-1))(h) \fr \Ho0C{3-s+h}\fr 0 $$\noindent
 where $h$ is the degree of $H$, 
we see that $H\eta\in\Hc0OQ{h-1}$. It follows that, if $Q'\in I_C$ is a surface without common components with $Q$, we have $H u (Q')=SQ' ({\rm mod}\ Q)$, where $S$ is a form of degree $h-1$;  since $H | Q$ we have $H | SQ'$; but $H$ and $Q'$ do not have common components and $\deg  
 S < \deg H$. It follows $S=0$, hence $H\eta =0$ and hence the image of $\eta$ is of the form 
 $\cJ_{E/Q}(-1)$, where $E\subset Q$ contains the surface $K$ where $H\cdot K=Q$. By [MDP] Proposition III.2.6 $\eta$
is not injective.\fine\par

\remk 1  Fix $s$ and $\xi\neq 0$. For  $Q\in\Hc0JCs$,  if the corresponding morphism $\eta: \cJ_{C/Q}\fr \cO_Q(-1))$  is  not injective we will denote by $H_Q$ the surface contained in $Q$  defined as the 2-dimensional component of the subscheme $D\subset Q$ whose ideal sheaf is $\ker \eta$;
we note that $H_Q$ can be characterized as follows:
$H_Q|Q$, $H_Q\cdot\xi=0$ and for every $H'$ s.t. $H'|Q$ and $H'\cdot\xi=0$ is $H_Q|H'$.
In fact from above result we have $H'\cdot\eta =0$ and hence the image of $\eta$ is of the form 
 $\cJ_{E/Q}(-1)$, where $E\subset Q$ contains the surface $K'$ whith $H'\cdot K'=Q$. Since $K$ is the 2-dimensional component of $E$ we have $K'|K$ and hence $H=H_Q$ divides $H'$.\par
 \proclaim{\bf Proposition 2}. Let $s(C)\leq s \leq e(C)+3$,
$Q,Q'\in\Hc0JCs$,  $\xi\in \Ho0C{3-s}$, $\xi\neq 0$ and let $\eta: \cJ_{C/Q}\fr \cO_Q(-1))$ and $\eta': \cJ_{C/Q'}\fr \cO_Q'(-1))$the morphisms  
 corresponding to $\xi$. Assume  $\eta$ and $\eta'$ are not injective and let $H_Q$ and 
$H_{Q'}$ the surfaces defined above. Then $H_Q=H_{Q'}$.\par
{\it Proof of Proposition 2}. Put $H=H_Q$ and $H'=H_{Q'}$. Assume they have degrees $0<h\leq h'<s$ .
Since $H\cdot\xi=0$ and $H'\cdot\xi=0$ by  remark 1 it is enough to prove $H|H'$.

Since $H'\cdot\xi=0$ we have $H'\cdot u (Q_i)=SQ_i ({\rm mod}\ Q)$ for all $Q_i\in I_C$, where $S$ is a surface of degree $h'-1$ and $u$ is as before; from this we get $K | S$ since the image $u(I_{C/Q})$ is contained in
the ideal $I_{E/Q}(-1)$ and hence all $u (Q_i)$ are multiple of $K$
; in particular for $Q_i=Q'$ , putting $u(Q)=KT$, $S=KV$ we have $H'T=H'K'V ({\rm mod}\ H)$ i.e. $H'T=H'K'V+ ZH$ with $Z$ of degree $h'-1$. It follows that $H'$ divides $ZH$ and hence $H'$ have a common factor of positive degree with $H$.\par Let $F=GCM(H,H')$;  
we will  show that $F=H$; assume the contrary, i.e. $\deg F <h$ and put
$H=F\cdot\overline{H}$, $Q=F\cdot\overline{Q}=F\cdot\overline{H}\cdot K$ and similarly
$H'=F\cdot\overline{H'}$, $Q'=F\cdot\overline{Q'}=F\cdot\overline{H'}\cdot K'$.\par\noindent
Now let $Y$ the residue curve of $C$ with respect to $F$ i. e. the curve whose homogeneous ideal is $(I_C:F)$: it is easy to see that $(I_C:F)$ is a saturated ideal and that 
 $Y\subset C$ is a curve i.e. does not have zero-dimensional components (isolated or embedded); moreover $Y$ is no empty: infact if it where $F\in I_C$ we would have $F\cdot\xi= 0$ and this is contrary to the fact that $H $ is the surface of minimal degree contained in $Q$ such that $H\cdot\xi= 0 $\par\noindent
From the exact sequence $$0\fr \cO_Y(-f)\fr \cO_C \fr \cO_{C\cap F}\fr 0$$
we have an exact sequence of sheaves on $C$
$$ 0\fr \omega_{C\cap F}\fr \omega_C \fr \omega_Y(f) $$
where $f=\deg F$.
From the above we obtain an exact sequence:
$$0\fr \Ho0{C\cap F}{3-s} \fr \Ho0{C}{3-s}\fr \Ho0{Y}{3-s+f}$$ and let $\overline{\xi}$ be the image of $\xi$ in $\Ho0{Y}{3-s+f}$.
Hence we have:\par\noindent
1) $\overline{\xi}\neq 0$:
In fact if it were zero, $\xi$ would be in $\Ho0{C\cap F}{3-s}$ and hence $F\cdot\xi= 0$; it follows $3-s+f\geq -e(Y)$\par\noindent
2) $Y\subset \overline{Q}$ and $Y\subset \overline{Q'}$:
it follows $s(Y)\leq s-f $.\par\noindent
3) $\overline{H}$ divides $\overline{Q}$ and has degree $0<h-f<s-f$;
$\overline{H'}$ divides $\overline{Q'}$ and has degree $0<h'-f<s-f$\par\noindent
4) $\overline{H}\cdot\overline{\xi}=0$ and $\overline{H'}\cdot\overline{\xi}=0$:
In fact $\Ho0{Y}{3-s+f}$ embeds in $\Ho0{C}{3-s+f}$ and the composition
 $\Ho0{C}{3-s}\fr \Ho0{Y}{3-s+f}\fr \Ho0{C}{3-s+f}$ is the moltiplication by $F$.
Hence $H\cdot\xi=\overline{H}\cdot F\xi=0$ implies that $\overline{H}\cdot\overline{\xi}=0$ and similarly for $H'$.\par\noindent
We have hence reproduced  for $Y$ the hypotheses of Proposition 2 and by the first part of proof $\overline{H}$ and $\overline{H'}$ have a common factor of positive degree and this is absurd since they were coprime.\fine\par

Now we conclude the proof that 2) implies 3) in  Theorem 1. 
Assume that  $C$ does not admit any elementary descending biliaison of height -1 and $s(C)\leq e(C)+3 $.  Using [MDP] Proposition III.2.6 three cases are possible.\par
i)  There is an $s$,  $s(C)\leq s\leq e(C)+3 $, a surface  $Q$  of degree $s$ containing $C$ and a surjective homomorphism $\cJ_{C/Q}\fr \cO_Q(-1)$; in this case $C$ is the section of $Q$ with a plane $S$ and hence we can make a descending elementary biliaison of $C$ on $S$ with a line.\par
ii) There is an $s$, $s(C)\leq s\leq e(C)+3 $, a surface  $Q$  of degree $s$ containing $C$ and  an injective non surjective homomorphism $\cJ_{C/Q}\fr \cO_Q(-1)$; in this case $C$ admit an elementary descendiong biliaison on $Q$.\par
iii) for all $s$, $s(C)\leq s\leq e(C)+3$ and all $Q$ of degree $s$, every non zero homomorphism $\cJ_{C/Q}\fr \cO_Q(-1)$ is not injective.
In this case, by Proposition 2, 
for every $s$ and every non zero $\xi\in\Ho0C{3-s}$ there is a form of degree $<s$ which divides every 
element in $\Hc0JCs$ and such that $H\xi =0$.\par
The last step is to make $H$ independent of $\xi$. Since for fixed $H$ the set of $\xi\in\Ho0C{3-s}$ s. t. $H\xi=0$ is a $\bf k$-subspace of $\Ho0C{3-s}$ and since the set of forms of degree $<s$ dividing the GCD of $\Hc0JCs$ is finite, we have the result.
This will finish our proof.\par

 \remk 2 We note that $H$ divides all  $\Hc0JC{i}$ for all $i=s(C), \dots, s$ and  kills all $\Ho0{C}{j}$ for $j=-e(C),\dots, 3-s$. This last statement follows from the fact that there is a linear form in $R$ which is not a zero divisor for $H_{*}(\omega_C)$.\par
\remk 3 If $C$ is a curve satisfying condition 3) c) of Theorem 1, then $C$ is not ACM. In fact an ACM curve has all its minimal generators in degree $\leq e(C)+3$ but for the condition 3) c) all these generators  have a common factor of positive degree. As a consequence we have that an ACM curve of degree $>1$ always admits a descending elementary biliaison.\par

\par\bigskip
\noindent
{\bf 2.} 
In this Section we prove the implication ${\rm 3)}\Rightarrow{\rm 1)}$ of Theorem 1.
Let $C$ be a curve satisfying condition 3) c); we want to prove that it is minimal in its biliaison class. \par
Let $$0\fr P\fr N_1\fr I_C\fr 0$$ be a minimal $N$-resolution of $I_C$ (see [LR] or [MDP]); we know that $P$ is free graded $R$-module. 
Now let $C'$ another curve in the same biliason class of $C$; by adding a free graded $R$-module $L$ we can assume that $C,C'$ have $N$-resolutions, (not necessarily minimal) of the form:
$$0\fr A\mapright{\alpha} N\mapright{\tau} I_C\fr 0$$
$$0\fr B\mapright{\beta} N\mapright{\sigma} I_{C'}(d)\fr 0$$
where
$$A=\bigoplus_{i=1}^n R(-a_i)=P\oplus L,\ a_1\leq a_2\leq\cdots\leq a_n$$
$$B=\bigoplus_{i=1}^n R(-b_i),\ b_1\leq b_2\leq\cdots\leq b_n$$
as in [LR] it is enough to prove that for all $i=1\cdots n$ is $a_i\leq b_i$; as in [LR] it is enough to prove that, for all $t\in \bZ$ is ${\rm rank}A_{\geq t}\geq {\rm rank}B_{\geq t}$; we consider three cases:
{\item 1\kern -5pt)}  $t\leq -(e(C)+4)$ in this case we use the same proof as in [LR];
{\item 2\kern -5pt)}  $t>-s(C)$ in this case we use the same proof as in [LR];
{\item 3\kern -5pt)}  let $-s(C)\geq t >-(e(C)+4)$ and assume that there exists $t$ such that ${\rm rank}A_{\geq t}< {\rm rank}B_{\geq t}$.
 We put $s=-t$ and we have $s(C)\leq s\leq e(C)+3$ and ${\rm rank}A_{\geq -s}< {\rm rank}B_{\geq -s}$. We assume also that $s$ is minimal with this property. 
\par\noindent
Put $B=B_1\oplus B_2$ where $B_1=\bigoplus_{b_i\leq s} B(-b_i)$,$B_2=\bigoplus_{b_i> s} B(-b_i)$
with ranks $r ,n-r$ and similarly $A=A_1\bigoplus A_2$ where $A_1=\bigoplus_{a_i\leq s} A(-a_i)$,$A_2=\bigoplus_{a_i> s} A(-a_i)$
with ranks $r' ,n-r'$ and $r>r'$.\par 

Let $H=H_s$ the form of degree $h<s$ given in 3)c) of Theorem 1; put $A'_2= A_2(h)$, 
$A'=A_1\oplus A'_2=A_1\oplus A_2(h)$, and define $\phi : A\fr A'$ by $\phi :A_1\oplus
A_2 \fr  A_1\oplus
A_2(h)$, $\phi =1\oplus  H$.\par 
 Having fixed notations we prove some Propositions.\par

\proclaim{\bf Proposition 3}. Under the above hypotheses and notations there is an exact sequence $0\bfr A_1\oplus
A'_2 \mapright{\alpha'} A'_2
\oplus M\mapright{\tau'} I_{C}\bfr 0 $ and a commutative diagram
$$\left.\matrix{0\fr &A&\hskip5pt\mapright{\alpha} N&\mapright{\tau} I_{C}&\fr 0\cr
&\hskip10pt\downarrow\phi&\hskip30pt\downarrow\psi&\hskip10pt\|&\cr
0\fr& A_1\oplus A'_2&\mapright{\alpha'} A'_2\oplus M&\mapright{\tau'} I_{C}&\fr 0\cr
}\right.$$
with $M$ torsion free graded $R$-module of rank $r'+1$ and where the map $\alpha'$
 is the direct sum of the identity on $A'_2$ and an inclusion $\theta :A_1\fr M$.\par

{\it Proof of Proposition 3}. The Proposition follows from the fact that, if $s<p$,
then $\Ext^1(I_C,R(-p))^0 \cong \Ho0C{4-p}$ with $3-s\geq 4-p$ is killed by $H$.
More in detail let $E\in \Ext^1(I_C,A)^0$ be the extension  $0\fr A\fr N\fr I_C\fr 0$ and let $\phi E\in\Ext^1(I_C,A')^0$ the composite extension (see [McL], Ch.III, Lemma 1.4); $\phi E$ has the form 
$$0\fr A_1\oplus A_2(h)\fr N'\fr I_C \fr 0 .$$
In order to prove our Proposition we prove that there is a map $N'\fr A_2(h)$ such that the composition $A_2(h)\fr N'\fr A_2(h)$ is the identity. In fact let $\overline{E}\in \Ext^1(I_C,A_2)^0$ the extension $0\fr A_2\fr \overline{N}\fr I_C\fr 0$ induced by $E$
, where $\overline{N}$ is the quotient of $N$ by $A_1 $; 
we see easily that the extension $H\cdot\overline{E}$ has the form: $$0\fr A_2(h)\fr N'/A_1\fr I_C\fr 0$$ and by the above observation it
splits; hence there is a map $N'/A_1\fr A_2(h)$ and the composition 
$A_2(h)\fr N' \fr N'/A_1\fr A_2(h)$ is the identity.  Moreover $A_1$ is contained in the kernel $M$ of the map $N'\fr A_2(h)$. \fine\par
By Proposition 3 we have hence the following exact sequences and commutative diagrams: 
$$\left.\matrix{0\fr &A&\mapright{\alpha} N&\mapright{\tau} I_{C}&\fr 0\cr
&\hskip10pt\downarrow\phi&\hskip20pt\downarrow\psi&\hskip10pt\|&\cr
0\fr& A'&\mapright{\alpha'} N'&\mapright{\tau'} I_{C}&\fr 0\cr
}\right.$$ where $N'=A'_2\oplus M$,

$$\left.\matrix{0\fr &A_2&\mapright{\alpha_2} N&\mapright{\rho} M&\fr 0\cr
&\hskip10pt\downarrow\cdot{\scriptstyle H}&\hskip20pt\downarrow\psi&\hskip10pt\|&\cr
0\fr& A'_2&\mapright{\alpha'_2} N'&\mapright{\rho'} M&\fr 0\cr
}\right.$$
where $\rho'$ is the projection, $\rho=\rho'\circ\psi$, $\alpha_2,\alpha'_2$ the restrictions of $\alpha,\alpha'$ respectively.\par

Moreover we have an exact sequence
 $$0\fr A_1\mapright{\theta} M\mapright{\lambda} I_C\fr 0$$ where $\lambda=\tau'_{\ \textstyle |M}$.\par
 \proclaim{\bf Proposition 4}. Under the above hypotheses and notations the composition
$$B_1\mapright{\beta_1} N\mapright{\psi} N' \mapright{\pi} A_2(h)$$ is zero, where 
$\beta_1=\beta_{\ \textstyle |B_1}$ and $\pi$ is the projection.
\par

{\it Proof of Proposition 4}. 
A summand of $B_1$ is of the form $R(-b)$ with $b\leq s$;  since the terms in $A_2(h)$
 have the form $R(-p+h)$, $s-h<p-h$, hence $-p+h+b<h$, the map $R(-b)\fr R(-p+h)$ will be zero if we will prove that it factors through the moltiplication by $H$ or equivalently that restricted to the surface $H$ it is zero.
To this end we first restrict to the affine open set $U={\bf P}^3\backslash S$ where $S$ is a surface containing $C$ and without common components with $H$. Since $\cJ_{C,U}\cong\cO_U$
the induced sequences $$\left.\matrix{0\bfr &\cA_U&\bfr\cN_U&\bfr\cJ_{C,U}&\fr 0\cr
&\hskip10pt\downarrow\tilde{\phi}_U&\hskip20pt\downarrow\tilde{\psi}_U&\|&\cr
0\fr&\cA_{1,U}\oplus \cA_2(h)_U&\fr \cA_2(h)_U\oplus\cM_U&\fr \cJ_{C,U}&\fr 0\cr
}\right.$$
split, hence they remain exact when restricted to $H\cap U$.
Now the map $\cO(-b)_U\fr \cJ_{C,U}$ induced by $R(-b)\fr I_C$ is zero when restricted to $H\cap U$, since $b\leq s$ and $H$ divides all $\Hc0JCb$ for $b\leq s$; so the map $\cO(-b)_U\fr \cO(-p+h)_U$  factors ( mod $H\cap U$) through  $\cA_U$, but the map
$\cA_U\fr \cO(-p+h)_U$ is zero ( mod $H\cap U$) .\fine\par
From the Proposition 4 follows that the map $\phi\circ\beta_1:B_1\fr N'$ factors though the inclusion $M\fr N'$ and since ${\rm rk }B_1=r>r'$ and
${\rm rk }M=r'+1$ it follows that ${\rm rk }B_1={\rm rk }M=r=r'+1$.
 \proclaim{\bf Proposition 5}. Under the above hypotheses and notations the composition map $$N \mapright{\sigma} I_{C'}(d) \mapright{\cdot H}I_{C'}(d+h)$$ factors through the map $\psi :N\fr N'$.\par
{\it Proof of Proposition 5}. We see easily that ${\rm coker }\psi$ is isomorphic to ${\rm coker }\phi$ hence it is annihilated by $H$. From the exact sequence $$ \Hom (N',I_{C'}(d))\fr \Hom (N,I_{C'}(d))\fr \Ext^1( {\rm coker }\psi,I_{C'}(d))$$ we get the result since $H$ annihilates $\Ext^1( {\rm coker }\psi,I_{C'}(d))$.\fine\par
Now let $\sigma' :N'\fr I_{C'}(d+h)$ be the map given by Proposition 5.
\proclaim{\bf Proposition 6}. Under the above hypotheses and notations the restriction 
 $\sigma'_{\ \textstyle |M} :M\fr I_{C'}(d+h)$ is zero.\par
{\it Proof of Proposition 6}. In fact the kernel contains $B_1$ hence has rank $r={\rm rk }M$ and $I_{C'}(d+h)$ is torsion free.\fine\par
\proclaim{\bf Proposition 7 }. Under the above hypotheses and notations the composition map 
 $\gamma: A_1\mapright{\alpha_1} N \mapright{\sigma} I_{C'}(d)$ is zero.\par
{\it Proof of Proposition 7 }. In fact the composition of $\gamma$ with the multiplication map 
$I_{C'}(d) \mapright{\cdot H}I_{C'}(d+h)$ factors through the map $\sigma'_{\ \textstyle |M} :M\fr I_{C'}(d+h)$ hence it is zero and also $\gamma$ is zero.\fine\par
As a conclusion the map $\alpha_1: A_1\fr N$ factors as $\beta_1\circ \zeta$ where $\zeta:A_1\fr B_1$ and $\beta_1:B_1\fr N$; it follows ${\rm rk}A_1\leq {\rm rk}B_1$. On the other hand remember that $s$ was minimum such that ${\rm rk}A_{\geq -s}< {\rm rk}B_{\geq -s}$. Hence we have:\par
a) for $t<s$ the map $\zeta$ induces an isomorphism $A_{\geq -t}\cong B_{\geq -t}$\par
b) for $t=s$ we can change  basis in $\bigoplus _{b_i=s}B(-b_i)$ such that $\bigoplus _{b_i=s}B(-b_i)=\bigoplus _{a_i=s}A(-a_i)\oplus R(-s)$ and $\zeta$ is the identity on $\bigoplus _{a_i=s}A(-a_i)$.\par
We can factor out $N$ by $A_1$ and we can assume in the sequel $A_1=0$ , $B_1=R(-s)$ $M=I_C$.\par

Let now $Q=H\cdot K$ be the image of the generator of $R(-s)$ in $I_C$; as in the proof of Proposition 4  we  restrict to the affine open set $U={\bf P}^3\backslash S$ where $S$ is a surface containing $C$ and without common components with $Q$. 
Moreover we still denote by $Q,H,K$ the local equations of the surfaces $Q,H,K$ on $U$.\par
Since $\cJ_{C,U}\cong\cO_U$
the induced sequences $0\fr \cA_U\fr\cN_U\fr\cJ_{C,U}\fr 0$ and $0\fr \cA'_U\fr\cN'_U\fr\cJ_{C,U}\fr 0$ 
split. Hence $\cN_U\cong \cN'_U\cong \cO_U^{\oplus n+1}$. We write down the matrix of the map $\psi_U:\cO_U^{\oplus n+1}\fr \cO_U^{\oplus n+1}$ induced by $\psi$. It has the form
$$\bf M=
\bordermatrix{&\grf{1}{A_2}&\grf{1}{I}\cr
{\scriptstyle A_2(h)\biggr \{}&{H\rm I_n}&{\rm C_n}\cr
{\hskip16pt\scriptstyle I \biggr\{}&0&1\cr}
$$
where $\rm I_n$ is the unit matrix and $\rm C_n$ is a one column matrix.\par
We also consider the matrix induced by the map $\beta: B\fr N$; it has the form 
$$\bf P=
\bordermatrix{&\grf{1}{B_1}&\grf{1}{B_2}\cr
{\hskip10pt\scriptstyle A_2\biggr \{}&{\rm D_n}&*\cr
{\hskip16pt\scriptstyle I \biggr\{}&Q&*\cr}
$$
where  $\rm D_n$  is a one column matrix.\par
The product $\bf M\cdot P$ is associated to the map $\tilde{\phi}_U\circ\tilde{\beta}_U : \cB_U\fr \cN'_U\cong \cA_2(h)_U\oplus \cJ_U$ and has the form
$$\bf M\cdot P=
\bordermatrix{&\grf{1}{B_1}&\grf{1}{B_2}\cr
{\hskip10pt\scriptstyle A_2(h)\biggr \{}&H{\rm D_n}+Q{\rm C_n}&*\cr
{\hskip26pt\scriptstyle I \biggr\{}&Q&*\cr}
$$
By Proposition 4 we have $H{\rm D_n}+Q{\rm C_n}=0$, hence
${\rm D_n}=-K{\rm C_n}$ and the matrix $\bf P$ has the form
$$\bf P=
\bordermatrix{&\grf{1}{B_1}&\grf{1}{B_2}\cr
{\hskip10pt\scriptstyle A_2\biggr \{}&-K{\rm C_n}&*\cr
{\hskip16pt\scriptstyle I \biggr\{}&Q&*\cr}
$$
We finish by observing that all $n\times n$ minors of $\bf P$ are multiple of $K$ and this is contrary to 
Hilbert-Burch Theorem.
\par\bigskip
\noindent


{\bf 3.}
In this section we deduce some consequences from Theorem 1. The first consequence,
whose proof is trivial, is an improvement of the Lazarsfeld-Rao property.
\proclaim{\bf Theorem 2}. Every curve in ${\bf P}^3$ can be obtained from a minimal one by means of a finite number of ascending elementary biliaisons.\par 
Next we recall that in  [MR] Maggioni and Ragusa determined bounds for the Betti numbers of a curve $C\subset{\bf P}^3$ that can be obtained as in Theorem 2,  so the results of [MR] apply to every curve in ${\bf P}^3$. In particular we have the following.
\proclaim{\bf Theorem 3}. Let $C$ any curve in ${\bf P}^3$ and let $C_0$ be a minimal curve in the biliaison class of $C$. Let us denote by $\nu(C), \nu(C_0)$ the minimal number of generators of $I_C, I_{C_0}$ respectively. Then we have $$\nu(C)\leq \nu(C_0)+ s(C)-s(C_0).$$\par
{\it Proof.} The proof is very easy; if $C$ is obtained from $C'$ by means of an elementary biliaison of type $(1,s)$ on a surface $Q$ of degree $s$, two cases are possible:\par
a) $s=s(C')$; in this case $s(C)=s(C')$ and $C,C'$ have the same number of minimal generators;\par
b) $s>s(C')$; in this case $s(C)=s(C')+1$ and $C$ has at most one more minimal generator than $C'$, i.e. the surface $Q$.\fine\par
We apply Theorem 3 to the following cases:\par
1) $C$ is a  ACM curve. In this case $C_0$ is a line and we get the well known bound $\nu(C)\leq s(C)+1$.\par
2) $C$ is an Arithmetically Buchsbaum curve of diameter $ \mu$. In this case it is $s(C_0)=2\mu$ and $\nu(C_0)=3\mu +1$ (see e.g. [MDP] ) and we find the Amasaki bound
$\nu(C)\leq s(C)+\mu +1$.\par

\par\bigskip 

\vskip 0.5in
\centerline{\bf REFERENCES}
\nobreak\vskip 0.5in\nobreak
\noindent\bib\hbox to 35 pt { [BBM]\hfill}E.Ballico, G.Bolondi, 
J.C.Migliore. 
{\it The Lazarsfeld-Rao problem for liaison classes of two-codimensional subschemes of ${\bf P}^n$}. Amer. J. Math. {\bf 113} (1991), 117-128. \par\medskip
\noindent\bib\hbox to 35 pt { [LR]\hfill}R.Lazarsfeld, A.P.Rao. 
{\it Linkage of general curves of large degree}. Lecture Notes in Math.,
vol. 156, Springer, 1970.\par\medskip
\noindent\bib\hbox to 35 pt { [McL]\hfill}S. Mac Lane. 
{\it Homology}, Springer, 1967.\par\medskip
\noindent\bib\hbox to 35 pt { [MDP]\hfill}M.Martin-Deschamps, D.Perrin. 
{\it Sur la clas\-sifica\-tion des cour\-bes gauches}. Ast\'e\-ris\-que {\bf 184-185}, 1990.\par\medskip
\noindent\bib\hbox to 35 pt { [MR]\hfill}R.Maggioni, A.Ragusa. 
{\it Betti numbers of space curves bounded by Hilbert functions}. Le Matematiche {\bf 52} (1997), 217-232 .
\par\medskip
\noindent\bib\hbox to 35 pt { [R]\hfill}A.P.Rao. 
{\it Liaison among curves in ${\bf P}^3$}. Inventiones Math. {\bf 50} (1979), 205-217.\par\medskip
\noindent\bib\hbox to 35 pt { [S]\hfill}E. Schlesinger. 
{\it The spectrum of projective curves}. Ph.D. Thesis, U.C.Berkeley, 1996.\par\medskip
\vskip 0.5in
\noindent Author address:\par
Dipartimento di Matematica\par
Universit\`a di Catania\par
Viale A.Doria, 6\qquad  95125 Catania (Italy)\par\medskip
\noindent E-mail:\par
sstrano @ dmi.unict.it
\bye